# Combined-distance-based score function of cognitive fuzzy sets and its application in lung cancer pain evaluation

Lisheng Jiang[1*], Tianyu Zhang[1], Shiyu Yan[1], Ran Fang[2]

[1] College of Management Science, Chengdu University of Technology, Chengdu, 610059, China

[2] School of Architecture and Civil Engineering, Chengdu University, Chengdu, 610106, China

**Abstract**

In decision making, the cognitive fuzzy set (CFS) is a useful tool in expressing experts' complex assessments of alternatives. The distance of CFS, which plays an important role in decision analyses, is necessary when the CFS is applied in solving practical issues. However, as far as we know, the studies on the distance of CFS are few, and the current Minkowski distance of CFS ignores the hesitancy degree of CFS, which might cause errors. To fill the gap of the studies on the distance of CFS, because of the practicality of the Hausdorff distance, this paper proposes the improved cognitive fuzzy Minkowski (CF-IM) distance and the cognitive fuzzy Hausdorff (CF-H) distance to enrich the studies on the distance of CFS. It is found that the anti-perturbation ability of the CF-H distance is stronger than that of the CF-IM distance, but the information utilization of the CF-IM distance is higher than that of the CF-H distance. To balance the anti-perturbation ability and information utilization of the CF-IM distance and CF-H distance, the cognitive fuzzy combined (CF-C) distance is proposed by establishing the linear combination of the CF-IM distance and CF-H distance. Based on the CF-C distance, a combined-distanced-based score function of CFS is proposed to compare CFSs. The proposed score function is employed in lung cancer pain evaluation issues. The sensitivity and comparison analyses demonstrate the reliability and advantages of the proposed methods.

## 1. Introduction

Decision making tries to give helpful suggestions to experts by analyzing the available information. In this regard, the information expression model and the decision analyzing techniques are two key points of decision making. In expressing information, we can use accurate numbers to describe the characteristic information of alternatives that can be measured in physical scales. For instance, the height of a student can be depicted as 1.6 meters. When the features of alternatives are difficult to be measured in physical scales, such as the individual potential, the subjective evaluations of experts might be helpful. For instance, let the degree to which a student belongs to the high potential group is between 0 and 1. High score means high potential. A teacher might give the evaluation that the degree to which a student belongs to the high potential group is 0.6. In depicting this kind of evaluation, the fuzzy set (Zadeh, 1965) is useful by employing the membership degree to represent the degree to which alternatives belong to certain sets. For the example of individual potential, the membership degree of the student belonging to the high potential group is 0.6.

---

* *Corresponding author*. Email: lsjiang96@163.com (Lisheng Jiang).

When evaluating alternatives, experts might give assessments in terms of both good and bad aspects. For instance, experts might think the degree to which alternative $A$ is good is 0.5 and the degree to which alternative $A$ is not good is 0.3. The numbers 0.5 and 0.3 are the positive information (good aspect) and negative information (bad aspect) of alternative $A$. For the fuzzy set, it is difficult to describe both positive and negative information. To model the two-dimensional positive and negative information, the intuitionistic fuzzy set (Atanassov, 1986) was proposed by considering the membership degree and non-membership degree to represent the positive and negative information. Because of the complexity of decision-making environments and the cognition limitation of experts, experts might occasionally give the special information such as the membership degree 0.5 and non-membership degree 0.6. Given that the intuitionistic fuzzy set needs the sum of membership and non-membership degrees to be not larger than 1, it cannot be used to deal with the special information. In this regard, two models: Pythagorean fuzzy set (Yager & Abbasov, 2013) and generalized orthopair fuzzy set (Yager & Alajlan, 2017) were proposed to handle the special information; nevertheless, they might cause the information loss, and they did not research the reasons of formation of the special information.

By further investigations, Jiang and Liao (2020) found that the knowledge and cognition limitation of experts might cause the cognitive overlaps of experts to good and bad judgements of alternatives. The overlaps of good and bad judgements of alternatives might cause the joint part of membership and non-membership degrees, resulting in the sum of membership and non-membership degrees larger than 1. In this scenario, Jiang and Liao (2020) proposed the cognitive fuzzy set (CFS) by taking into account the joint part called the joint degrees of membership and non-membership degrees. Compared with the Pythagorean fuzzy set (Yager & Abbasov, 2013) and generalized orthopair fuzzy set (Yager & Alajlan, 2017), the CFS has high interpretability and low information loss (Jiang & Liao, 2020). Although the CFS has been studied in other theories like the evidential reasoning approach (Zeng *et al.*, 2022), the distance measure of CFS is not well studied.

For a mass of decision analyzing techniques such as the consistency repairing (Yang *et al.*, 2019; Meng *et al.*, 2020; Kumar & Chen, 2022; Lu *et al.*, 2023), consensus improvement (Wang *et al.*, 2018; Li *et al.*, 2020; Lu *et al.*, 2023; Yao & Xu, 2023; Zhang & Dai, 2023), decision-making methods (Chen, 2000; Opricovic & Tzeng, 2004; Gou *et al.*, 2020), and regression models (Grama & Neumann, 2006; Singh & Marx, 2013), the distance is one of the foundations. For the distance of CFS, Jiang and Liao (2020) proposed the Minkowski distance of CFS, but the hesitancy degree of CFS was not involved in the distance, which might lead to errors. Besides the Minkowski distance of CFS (Jiang & Liao, 2020), as far as we know, there are few studies on the distance of CFS, which might limit the applications and developments of CFS.

Given that the studies on the distance of CFS are important but not adequate, because of the practicality of the Hausdorff distance (Yang & Ding, 2023; Zhang *et al.*, 2023; Karaaslan & Karamaz, 2024; Liu *et al.*, 2025), this paper aims to research the Hausdorff distance of CFS and improve the Minkowski distance of CFS proposed by Jiang and Liao (2020). To achieve this goal, we propose the cognitive fuzzy Hausdorff (CF-H) distance based on the interval-expressed CFSs and the improved cognitive fuzzy Minkowski (CF-IM) distance by considering the hesitancy degree of CFS. Then, the cognitive fuzzy combined (CF-C) distance is introduced to balance the anti-perturbation ability and the information utilization of the CF-H distance and CF-IM distance. A score function considering the CF-C distance is given to compare CFSs.

The innovations and contributions of this paper can be summarized as follows:

(1) The CF-H distance and CF-IM distance are studied. The hesitancy degree of CFS is considered to improve the Minkowski distance proposed by Jiang and Liao (2020). Then, the CF-H distance is researched based on the interval-expressed CFS. A simulation is done to show the characteristics of the anti-perturbation ability and information utilization of the CF-H distance and CF-IM distance.

(2) The combined-distance-based score function of CFS is proposed. To balance the anti-perturbation ability and information utilization of the CF-H distance and CF-IM distance, the CF-C distance is proposed by employing the linear combination of the CF-H distance and CF-IM distance. Based on the CF-C distance, the combined-distance-based score function of CFS is proposed. Three examples are given to show the relations and features of the CF-H distance, CF-IM distance, and CF-C distance.

(3) An illustration on the lung cancer pain evaluation gives an insight into the applicability of the combined-distance-based score function. The sensitivity analysis demonstrates the reliability of the proposed method. The comparison analysis is done to show the advantages of the proposed method.

This paper is organized as follows: Section 2 reviews the relevant studies on the CFS and the Hausdorff distance of intuitionistic fuzzy sets. Section 3 proposes the CF-H distance, CF-IM distance, and CF-C distance. The combined-distance-based score function is also defined in Section 3. The case study is done in Section 4. This paper closes with conclusions in Section 5.

## 2. Preliminaries

This section reviews the relevant knowledge about the CFS and the Hausdorff distance of intuitionistic fuzzy sets to facilitate the further discussions.

### 2.1. The concepts of CFSs

To depict the twofold positive and negative information in evaluation, Atanassov (1986) proposed the intuitionistic fuzzy set consisting of membership degree $u$ (positive information) and non-membership degree $v$ (negative information) with the restriction $u+v \leq 1$. The hesitancy degree is defined as $h=1-u-v$. The intuitionistic fuzzy set has been developed in theory (Ouyang & Pedrycz, 2016; Zhou et al., 2020) and has good performance in practical issues such as the supplier selection (Wan & Li, 2013) and company evaluation (Zhou & Xu, 2020). With the applications of intuitionistic fuzzy set, researchers found that when experts give the information like $u=0.5$ and $v=0.6$, the intuitionistic fuzzy set cannot be employed. To solve this issue, the Pythagorean fuzzy set (Yager & Abbasov, 2013) and generalized orthopair fuzzy set were proposed (Yager & Alajlan, 2017). Although the pythagorean fuzzy set and generalized orthopair fuzzy set worked well in mathematic operations, the parameters of them had low interpretability (Jiang & Liao, 2020), and there might be the information loss when doing computations (Yager & Alajlan, 2017).

To overcome the above drawbacks, the CFS was proposed by considering the joint part of $u$ and $v$ (Jiang & Liao, 2020). Specifically, because of the complexity of decision-making environments and the cognitive limitations of human, experts might be difficult to clearly distinguish the positive and negative information, resulting in the joint part of membership and non-membership degrees. When the joint part was large, the sum of the membership and non-membership degree might be larger than 1 (Jiang & Liao, 2020). In this regard, The joint degree $j$ of $u$ and $v$ was defined to represent the size of joint part, which can be explained as the confused level of experts when giving information (Jiang & Liao, 2020).

Let the set of evaluated alternatives be $T = \{t_1, t_2, \cdots, t_n\}$. The CFS is defined as $F = \{t, <u(t), v(t), j(t)> | t \in T\}$ where $u(t)$, $v(t)$, and $j(t)$ belong to [0,1]. $j(t)$ is the joint degree of $u(t)$ and $v(t)$. The true membership and non-membership degrees are $u^*(t) = u(t) - j(t)$ and $v^*(t) = v(t) - j(t)$, satisfying the inequalities $0 \leq u^*(t) + v^*(t) \leq 1$. A CFS can be abbreviated as $f_t = <u_t, v_t, j_t>$ called the cognitive fuzzy number (CFN). It is easy to get $j_t \in [\max\{0, u_t + v_t - 1\}, \min\{u_t, v_t\}]$. The hesitancy degree $h_t$ defined as $1 - u_t - v_t + j_t$ describes the degree to which experts have insufficient information about the alternative $t$. The relations of $u_t$, $v_t$, $j_t$, and $h_t$ are intuitively shown in Figure 1.

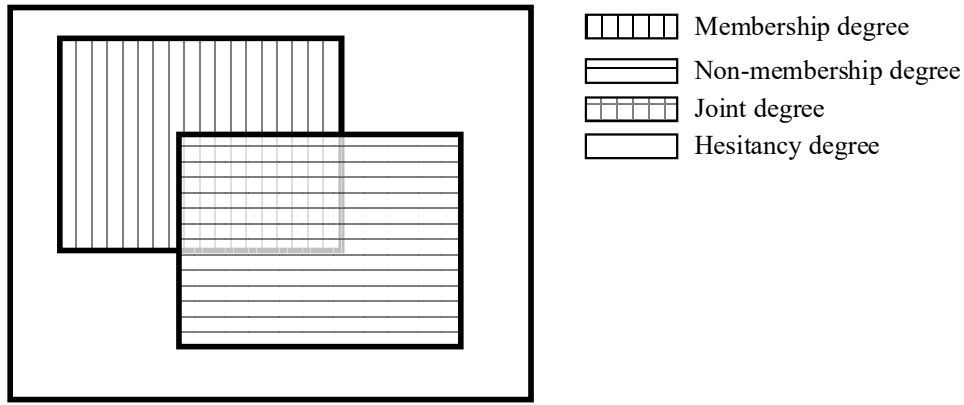

**Figure 1.** The relations of $u_t$, $v_t$, $j_t$, and $h_t$.

To apply the CFS in real issues, Jiang and Liao (2020) proposed the Minkowski distance of two CFNs $f_1 = <u_1, v_1, j_1>$ and $f_2 = <u_2, v_2, j_2>$ as Eq. (1), where $p$ is the positive integer. When $p = 1$, $|u_1^* - u_2^*| + |v_1^* - v_2^*| + |j_1 - j_2|$ is the Manhattan distance. When $p = 2$, $\sqrt{|u_1^* - u_2^*|^2 + |v_1^* - v_2^*|^2 + |j_1 - j_2|^2}$ is the Euclidean distance. When $p \to +\infty$, $\max(|u_1^* - u_2^*|, |v_1^* - v_2^*|, |j_1 - j_2|)$ is the Chebyshev distance.

$$d_m'(f_1, f_2) = (|u_1^* - u_2^*|^p + |v_1^* - v_2^*|^p + |j_1 - j_2|^p)^{1/p} \qquad (1)$$

**Example 1.** For two CFNs $f_1 = <0.3, 0.2, 0.1>$ and $f_2 = <0.4, 0.3, 0.1>$, according to Eq. (1), the distances of them to the CFN $f_3 = <1, 0, 0>$ are both 1 when $p = 1$, so the two CFNs $f_1$ and $f_2$ have the same distance to $f_3$. If we consider the hesitancy degree, the hesitancy degree of $f_2$ shown as $h_2 = 0.4$ is closer to the hesitancy degree of $f_3$ shown as $h_3 = 0$ than the hesitancy degree of $f_1$ shown as $h_1 = 0.6$. Hence, $f_2$ is closer to $f_3$ than $f_1$. This example shows that the ignorance of hesitancy degree might cause inaccurate results.

Because the joint degree in CFS has the practical meaning, which is the confused level of experts, the CFS has higher interpretability than the Pythagorean fuzzy set and generalized orthopair fuzzy set. Although the CFS has the above advantages, the Minkowski distance of CFS needs further investigations. In Eq. (1), the hesitancy degree of CFS is not considered, resulting in the inaccuracy of the distance of CFS. The inaccurate distance might limit the application of CFS.

## 2.2. The Hausdorff distance of intuitionistic fuzzy sets

In the intuitionistic fuzzy set, $u$ is the degree of one alternative being good. $h$ is the hesitancy degree meaning that the alternative might be good or not good. By considering $h$, the possible maximum degree of one alternative being good is $u + h = 1 - v$. In this regard, the intuitionistic fuzzy set is represented by an interval $[u, 1-v]$ (Grzegorzewski, 2004; Xu, 2007).

The general expression of the Hausdorff distance between two sets $\Lambda_1$ and $\Lambda_2$ is defined as Eq. (2) (Grzegorzewski, 2004). Eq. (2) tries to find the maximum distance from one point in $\Lambda_1$ to the nearest point in $\Lambda_2$. When the two sets are interval-expressed intuitionistic fuzzy sets $[u_1, 1-v_1]$ and $[u_2, 1-v_2]$, the Hausdorff distance of intuitionistic fuzzy set is defined as Eq. (3) (Grzegorzewski, 2004), where $u_1$ and $v_1$ are the membership and non-membership degrees of the first intuitionistic fuzzy set, and $u_2$ and $v_2$ are the membership and non-membership degrees of the second intuitionistic fuzzy set.

$$d_h'(\Lambda_1, \Lambda_2) = \max\{\sup_{\varsigma \in \Lambda_1} \inf_{\xi \in \Lambda_2} |\varsigma - \xi|, \sup_{\xi \in \Lambda_2} \inf_{\varsigma \in \Lambda_1} |\varsigma - \xi|\} \qquad (2)$$

$$d_h' = \max\{|u_1 - u_2|, |v_1 - v_2|\} \qquad (3)$$

The Hausdorff distance has been widely applied in real areas such as the multi-criteria decision making (Li *et al.*, 2015), medical diagnosis (Luo *et al.*, 2020), and hospital management (Liao *et al.*, 2019). These applications demonstrated the applicability and reliability of the Hausdorff distance.

## 3. The combined-distance-based score function of CFS

Because the Hausdorff distance is a widely used distance, this section proposes the CF-H distance. The CF-IM distance is introduced by adding the hesitancy degree of CFS to the Minkowski distance proposed by Jiang and Liao (2020). It is found that the CF-H distance has strong anti-perturbation ability but low information utilization while the CF-IM distance has high information utilization but weak anti-perturbation ability. As the CF-H distance and CF-IM distance can make up for each other's shortcomings, this section uses them to propose the CF-C distance based on which the combined-distance-based score function of CFN is proposed to compare CFNs.

### 3.1. The CF-C distance

*3.1.1. The CF-IM distance*

Let two CFNs be $f_{t_1} = <u_{t_1}, v_{t_1}, j_{t_1}>$ and $f_{t_2} = <u_{t_2}, v_{t_2}, j_{t_2}>$. By considering the hesitancy degree, the CF-IM distance can be defined as Eq. (4), where $p$ is the positive integer. $u^*_{t_1} = u_{t_1} - j_{t_1}$, $v^*_{t_1} = v_{t_1} - j_{t_1}$, $u^*_{t_2} = u_{t_2} - j_{t_2}$, $v^*_{t_2} = v_{t_2} - j_{t_2}$, $h_{t_1} = 1 - u_{t_1} - v_{t_1} + j_{t_1}$, and $h_{t_2} = 1 - u_{t_2} - v_{t_2} + j_{t_2}$.

$$d_m(f_{t_1}, f_{t_2}) = (|u^*_{t_1} - u^*_{t_2}|^p + |v^*_{t_1} - v^*_{t_2}|^p + |j_{t_1} - j_{t_2}|^p + |h_{t_1} - h_{t_2}|^p)^{1/p} \quad (4)$$

**Property 1.** (1) $d_m(f_{t_1}, f_{t_1}) = 0$; (2) $d_m(f_{t_1}, f_{t_2}) = d_m(f_{t_2}, f_{t_1})$; (3) $d_m(f_{t_1}, f_{t_2}) + d_m(f_{t_2}, f_{t_3}) \geq d_m(f_{t_1}, f_{t_3})$.

**Proof.** (1) and (2) are obviously established. For (3), let $\vec{f}_{t_1} = \overrightarrow{(u^*_{t_1}, v^*_{t_1}, j_{t_1}, h_{t_1})}$, $\vec{f}_{t_2} = \overrightarrow{(u^*_{t_2}, v^*_{t_2}, j_{t_2}, h_{t_2})}$, and $\vec{f}_{t_3} = \overrightarrow{(u^*_{t_3}, v^*_{t_3}, j_{t_3}, h_{t_3})}$. $d_m(f_{t_1}, f_{t_3}) = \|\vec{f}_{t_1} - \vec{f}_{t_3}\|_p = \|\vec{f}_{t_1} - \vec{f}_{t_2} - \vec{f}_{t_3} + \vec{f}_{t_2}\|_p \leq \|\vec{f}_{t_3} - \vec{f}_{t_2}\|_p + \|\vec{f}_{t_1} - \vec{f}_{t_2}\|_p = d_m(f_{t_1}, f_{t_2}) + d_m(f_{t_2}, f_{t_3})$. This completes the proof.

*3.1.2. The CF-H distance*

In the CFN $f_t = <u_t, v_t, j_t>$, $u^*_t$ represents the true degree of alternative $t$ being good. For $j_t$ and $h_t$, there is the possibility that alternative $t$ is good. By considering $j_t$ and $h_t$, the maximum degree of $t$ being good is $u_t^* + j_t + h_t = 1 - v_t^*$. Then, the CFN can be represented as an interval $[u_t^*, 1-v_t^*]$. Based on the interval-expressed CFN, the CF-H distance $f_{t_1} = <u_{t_1}, v_{t_1}, j_{t_1}>$ and $f_{t_2} = <u_{t_2}, v_{t_2}, j_{t_2}>$ is defined as Eq. (5), where $u_{t_1}^* = u_{t_1} - j_{t_1}$, $v_{t_1}^* = v_{t_1} - j_{t_1}$, $u_{t_2}^* = u_{t_2} - j_{t_2}$, $v_{t_2}^* = v_{t_2} - j_{t_2}$.

$$d_h(f_{t_1}, f_{t_2}) = \max\{|u_{t_1}^* - u_{t_2}^*|, |v_{t_1}^* - v_{t_2}^*|\} \tag{5}$$

**Property 2.** (1) $d_h(f_{t_1}, f_{t_1}) = 0$; (2) $d_h(f_{t_1}, f_{t_2}) = d_h(f_{t_2}, f_{t_1})$; (3) $d_h(f_{t_1}, f_{t_2}) + d_h(f_{t_2}, f_{t_3}) \geq d_h(f_{t_1}, f_{t_3})$.

**Proof.** (1) and (2) are obviously established. For (3), it is easy to get $\max\{|u_{t_1}^* - u_{t_2}^*|, |v_{t_1}^* - v_{t_2}^*|\} + \max\{|u_{t_3}^* - u_{t_2}^*|, |v_{t_3}^* - v_{t_2}^*|\} \geq |u_{t_1}^* - u_{t_2}^*| + |u_{t_3}^* - u_{t_2}^*|, |v_{t_1}^* - v_{t_2}^*| + |v_{t_3}^* - v_{t_2}^*|$. Because $|u_{t_1}^* - u_{t_2}^*| + |u_{t_3}^* - u_{t_2}^*| \geq |u_{t_1}^* - u_{t_3}^*|$ and $|v_{t_1}^* - v_{t_2}^*| + |v_{t_3}^* - v_{t_2}^*| \geq |v_{t_1}^* - v_{t_3}^*|$, $\max\{|u_{t_1}^* - u_{t_2}^*|, |v_{t_1}^* - v_{t_2}^*|\} + \max\{|u_{t_3}^* - u_{t_2}^*|, |v_{t_3}^* - v_{t_2}^*|\} \geq |u_{t_1}^* - u_{t_3}^*|, |v_{t_1}^* - v_{t_3}^*|$. Then, it can be inferred that $\max\{|u_{t_1}^* - u_{t_2}^*|, |v_{t_1}^* - v_{t_2}^*|\} + \max\{|u_{t_3}^* - u_{t_2}^*|, |v_{t_3}^* - v_{t_2}^*|\} \geq \max\{|u_{t_1}^* - u_{t_3}^*|, |v_{t_1}^* - v_{t_3}^*|\}$, so $d_h(f_{t_1}, f_{t_2}) + d_h(f_{t_2}, f_{t_3}) \geq d_h(f_{t_1}, f_{t_3})$. This completes the proof.

Because experts are confused when giving the CFS, there might be the errors in the given information. Let the errors be the perturbations on membership and non-membership degrees. An example is given to investigate the influence of the perturbations on the CF-H distance and CF-IM distance by a simulation experiment.

**Example 2.** Let two CFNs be $f_1 = <0.8, 0.4, 0.32>$ and $f_2 = <0.1, 0.9, 0.09>$. The CF-IM distances of them are 1.46, 0.8987, and 0.7964 for $p = 1, 2, 3$ calculated by Eq. (4). The CF-H distance of them is 0.73 computed by Eq. (5). Suppose there is a perturbation $\varepsilon$ on the membership and non-membership degrees of $f_1$. Then, $f_1$ is changed to $<0.8+\varepsilon, 0.4-\varepsilon, 0.32>$. Because $0.4-\varepsilon \geq 0.32$ and $0.8+\varepsilon \geq 0.32$, $\varepsilon \in [-0.48, 0.08]$. A 100-time simulation is done by randomly getting $\varepsilon$ from $[-0.48, 0.08]$ to uncover the perturbations on the CF-IM distance and CF-H distance. The results are shown in Figure 2.

In Figures 2(a) and 2(d), $\Delta d_m$ is the absolute values of the perturbations on the CF-IM distance with $p = 1$. $\Delta d_h$ is the absolute values of the perturbations on the CF-H distance. Because $\Delta d_m - \Delta d_h \geq 0$ known from Figure 2(d), it is clear that the CF-IM distance has higher perturbations than the CF-H distance. We can have the same conclusion by the vertical analyses to the figures in Figure 2. By the horizontal comparisons to the figures in Figure 2, it is found that when $p$ rises, the gap between $\Delta d_m$ and $\Delta d_h$ reduces.

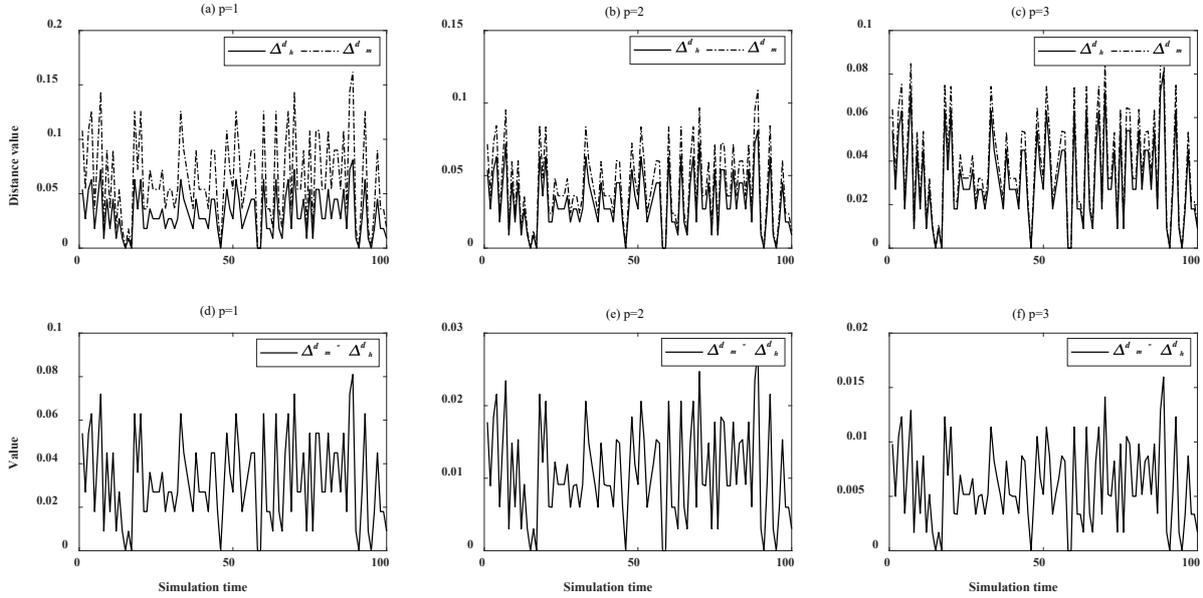

**Figure 2.** The perturbations on the CF-IM distances and CF-H distance.

**Insight 1.** The simulation in Example 2 demonstrates that the anti-perturbation ability of the CF-H distance is stronger than that of the CF-IM distance.

Although the anti-perturbation ability of the CF-H distance is strong, the maximum operation in Eq. (5) might lose part of the evaluation information in computing. On the contrary, all information, such as membership degrees, non-membership degrees, joint degrees, and hesitancy degrees, is reflected in Eq. (4), so the information utilization of the CF-IM distance is high. Given that the CF-H distance and the CF-IM distance might overcome each other's drawbacks, the CF-C distance of two CFNs $f_{t_1} = <u_{t_1}, v_{t_1}, j_{t_1}>$ and $f_{t_2} = <u_{t_2}, v_{t_2}, j_{t_2}>$ is given as Eq. (6), where $\lambda \in [0,1]$ is the balance parameter. $d_m(f_{t_1}, f_{t_2})$ and $d_h(f_{t_1}, f_{t_2})$ are the CF-IM distance and CF-H distance of $f_{t_1} = <u_{t_1}, v_{t_1}, j_{t_1}>$ and $f_{t_2} = <u_{t_2}, v_{t_2}, j_{t_2}>$.

$$d_c(f_{t_1}, f_{t_2}) = \lambda d_m(f_{t_1}, f_{t_2}) + (1-\lambda) d_h(f_{t_1}, f_{t_2}) \tag{6}$$

**Property 3.** (1) $d_c(f_{t_1}, f_{t_1}) = 0$; (2) $d_c(f_{t_1}, f_{t_2}) = d_c(f_{t_2}, f_{t_1})$; (3) $d_c(f_{t_1}, f_{t_2}) + d_c(f_{t_2}, f_{t_3}) \geq$

$d_c(f_{t_1}, f_{t_3})$.

**Proof.** According to Property 1 and Property 2, (1) and (2) are obviously established. For (3), $d_c(f_{t_1}, f_{t_2}) + d_c(f_{t_2}, f_{t_3}) = \lambda(d_m(f_{t_1}, f_{t_2}) + d_m(f_{t_3}, f_{t_2})) + (1-\lambda)(d_h(f_{t_1}, f_{t_2}) + d_h(f_{t_3}, f_{t_2}))$. According to Property 1 and Property 2, $\lambda(d_m(f_{t_1}, f_{t_2}) + d_m(f_{t_3}, f_{t_2})) + (1-\lambda)(d_h(f_{t_1}, f_{t_2}) + d_h(f_{t_3}, f_{t_2})) \geq \lambda d_m(f_{t_1}, f_{t_3}) + (1-\lambda)d_h(f_{t_1}, f_{t_3}) = d_c(f_{t_1}, f_{t_3})$. Hence, $d_c(f_{t_1}, f_{t_2}) + d_c(f_{t_2}, f_{t_3}) \geq d_c(f_{t_1}, f_{t_3})$. This completes the proof.

**Theorem 1.** When $p = 1$, $d_m(f_{t_1}, f_{t_2}) \geq d_c(f_{t_1}, f_{t_2}) \geq d_h(f_{t_1}, f_{t_2})$.

**Proof.** When $p = 1$, $d_m(f_{t_1}, f_{t_2}) = |u^*_{t_1} - u^*_{t_2}| + |v^*_{t_1} - v^*_{t_2}| + |j_{t_1} - j_{t_2}| + |h_{t_1} - h_{t_2}| \geq \max\{|u^*_{t_1} - u^*_{t_2}|, |v^*_{t_1} - v^*_{t_2}|\} = d_h(f_{t_1}, f_{t_2})$. $d_m(f_{t_1}, f_{t_2}) - d_c(f_{t_1}, f_{t_2}) = (1-\lambda)(d_m(f_{t_1}, f_{t_2}) - d_h(f_{t_1}, f_{t_2})) \geq 0$, so $d_m(f_{t_1}, f_{t_2}) \geq d_c(f_{t_1}, f_{t_2})$. $d_c(f_{t_1}, f_{t_2}) - d_h(f_{t_1}, f_{t_2}) = \lambda(d_m(f_{t_1}, f_{t_2}) - d_h(f_{t_1}, f_{t_2})) \geq 0$, so $d_c(f_{t_1}, f_{t_2}) \geq d_h(f_{t_1}, f_{t_2})$. In a word, $d_m(f_{t_1}, f_{t_2}) \geq d_c(f_{t_1}, f_{t_2}) \geq d_h(f_{t_1}, f_{t_2})$. This completes the proof.

**Example 3.** Let $\lambda = 0.5$. Using the same two CFNs in Example 1, the CF-C distance of them is 1.095. When $\lambda$ changes from 0 to 1, the trends of three types of distances are shown in Figure 3. Theorem 1 is verified from the results of the simulation in Figure 3(a).

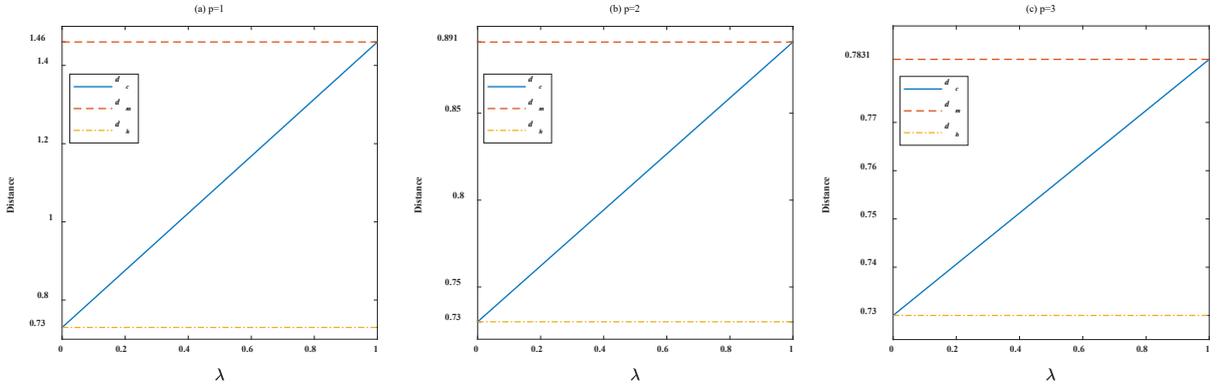

**Figure 3.** The trends of three types of distances when $\lambda$ changes from 0 to 1.

**Example 4.** Similar to Example 1, suppose there is a perturbation $\varepsilon \in [-0.48, 0.08]$ on the membership and non-membership degrees of $f_1$. Then, $f_1$ is changed to $<0.8+\varepsilon, 0.4-\varepsilon, 0.32>$. Another simulation is done by randomly getting $\varepsilon$ from $[-0.48, 0.08]$ to uncover the perturbations on the CF-IM distance, CF-H distance, and CF-C distance. The results are shown in Figure 4.

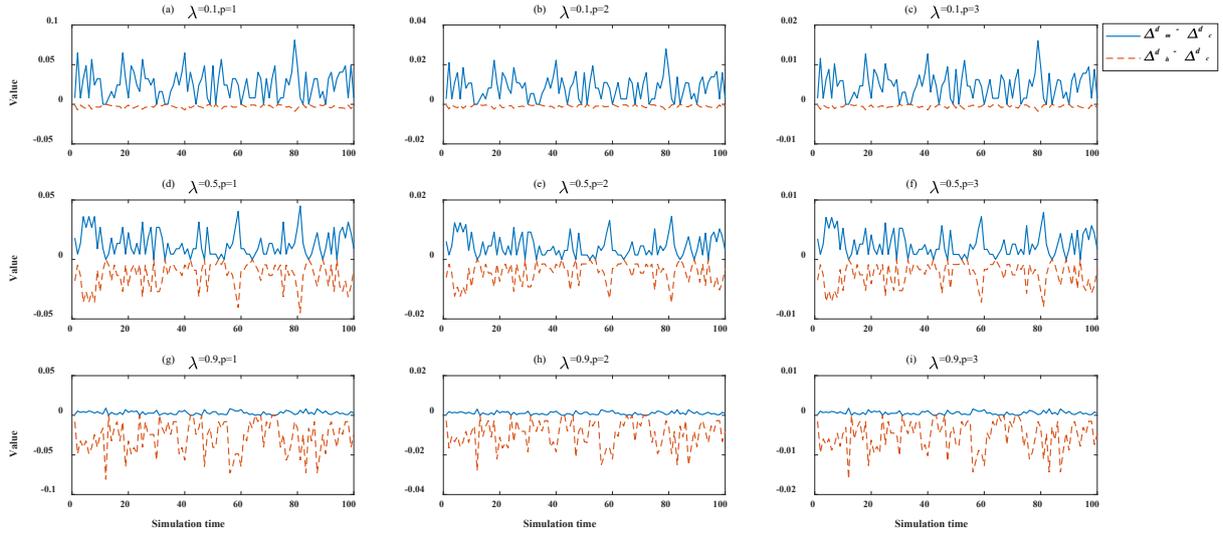

**Figure 4.** The relations of the perturbations on the three types of distances of different $\lambda$ and $p$.

In Figure 4, $\Delta d_h$ is the absolute values of the perturbations on the CF-H distance. $\Delta d_m$ is the absolute values of the perturbations on the CF-IM distance. $\Delta d_c$ is the absolute values of the perturbations on the CF-C distance. By the vertical comparisons to the figures in Figure 3, it is found that when $\lambda$ rises, $\Delta d_c$ approaches $\Delta d_m$. Because $\Delta d_m \geq \Delta d_h$ known from Example 1, $\Delta d_c$ becomes large when $\lambda$ rises, which shows that the anti-perturbation ability of the CF-C distance reduces with the increase of $\lambda$. It is found that the above conclusion about the anti-perturbation ability of the CF-C distance is correct when $p = 1, 2, 3$. Meanwhile, by the horizontal comparisons of figures in Figure 3, we find that $\Delta d_m - \Delta d_c$ and $\Delta d_c - \Delta d_h$ reduce when $p$ rises.

Examples 2 to 4 give two insights into the CF-IM distance, the CF-H distance, and the CF-C distance shown as follows.

**Insight 2.** The anti-perturbation ability of the CF-C distance descends with the rise of the balance parameter $\lambda$.

**Insight 3.** When $p$ increases, $\Delta d_m - \Delta d_c$, $\Delta d_c - \Delta d_h$, and $\Delta d_m - \Delta d_h$ reduce. The anti-perturbation abilities of the CF-IM distance, CF-H distance, and CF-C distance tend to be similar.

### 3.2. The combined-distance-based score function of CFN

For the convenience of comparing CFNs, this section proposes the combined-distance-based score function of CFNs.

For the CFN $<1,0,0>$, the membership degree is 1 without confusion and hesitancy, so $<1,0,0>$ is the CFN with best performance. On the contrary, the non-membership degree of $<0,1,0>$ is 1 without confusion and hesitancy, so $<0,1,0>$ is the CFN with worst performance. For any CFN $f_t$, the farther $f_t$ is from $<0,1,0>$, the better $f_t$ is. In this regard, motivated by the TOPSIS method (Zhang *et al.*, 2020; Singh & Gupta, 2020), the combined-distance-based score function of CFN $f_t =<u_t,v_t,j_t>$ is defined as Eq. (7), where $\bar{d}_c(f_t) = \lambda d_m(f_t,<0,1,0>) + (1-\lambda)d_h(f_t,<0,1,0>)$ is the CF-C distance between $f_t =<u_t,v_t,j_t>$ and $<0,1,0>$. $\tilde{d}_c(f_t) = \lambda d_m(f_t,<1,0,0>) + (1-\lambda)d_h(f_t,<1,0,0>)$ is the CF-C distance between $f_t =<u_t,v_t,j_t>$ and $<1,0,0>$. $\bar{d}_c(f_t) + \tilde{d}_c(f_t)$ is used to normalize the score. It is easy to get $s_t \in [0,1]$.

$$s_t = \frac{\bar{d}_c(f_t)}{\bar{d}_c(f_t) + \tilde{d}_c(f_t)} \tag{7}$$

High score means the CFN is close to $<1,0,0>$ having the high membership degree (positive information) and good performance. For two CFNs $f_1$ and $f_2$, if $s_1 > s_2$, $f_1$ is better than $f_2$. If $s_1 = s_2$, $f_1$ and $f_2$ have the same performance.

**Example 5.** When $\lambda = 0.5$ and $p = 2$, the scores of the two CFNs $f_1 =<0.8,0.4,0.32>$ and $f_2 =<0.1,0.9,0.09>$ in Example 1 are $s_1 = 0.6369$ and $s_2 = 0.1555$, so $f_1$ is better than $f_2$. The scores with different $\lambda$ and $p$ are shown in Figure 5.

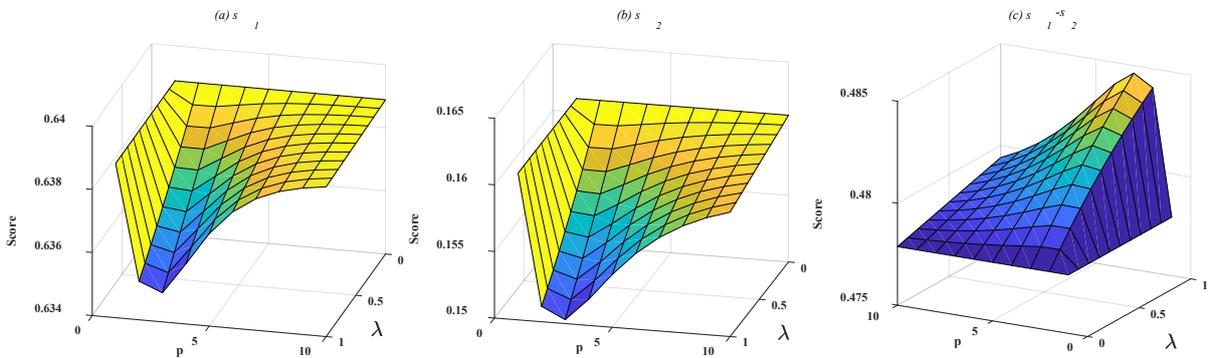

**Figure 5.** The scores with different $\lambda$ and $p$.

Figure 5(a) and Figure 5(b) show the scores of $f_1$ and $f_2$ when $\lambda \in [0,1]$ and $p=1,2,3,4,5,6,7,8,9,10$. Figure 5(c) shows that when $\lambda$ and $p$ change, $s_1$ is always larger than $s_2$, so the conclusion that $f_1$ is better than $f_2$ is reliable. Example 4 gives the insight into the high robustness of the comparison method based on the combined-distance-based score function.

## 4. Case study on the lung cancer pain evaluation

This section uses the combined-distance-based score function of CFN to solve the lung cancer pain evaluation issues, so as to show the applicability of the proposed method. The sensitivity and comparison analyses are done to illustrate the reliability of the results.

### 4.1. Problem description

After the lung resection or in the advanced lung cancer, the pain caused by the cancer has obvious interference on the recovery and life quality of patients, resulting in the burden on the patients' body and mind. In this regard, the National Comprehensive Cancer Network (2020) suggested that patients should be evaluated for pain, and then the treatment program should be improved according to different pain levels to improve the life quality of patients.

On the one hand, the pain evaluation is subjective, so patients' feel must be considered. The National Comprehensive Cancer Network (2020) gives a list shown in Table 1 to help patients evaluate pain.

**Table 1.** The impact of pain evaluation.

| 1. General Activity |
|---|
| 0   1  2  3  4  5  6  7  8  9   10 |
| Does not interfere                       Completely interfere |
| 2. Mood |
| 0   1  2  3  4  5  6  7  8  9   10 |
| Does not interfere                       Completely interfere |
| 3. Walking Ability |
| 0   1  2  3  4  5  6  7  8  9   10 |
| Does not interfere                       Completely interfere |
| 4. Normal Work (including both work outside the home and housework) |
| 0   1  2  3  4  5  6  7  8  9   10 |
| Does not interfere                       Completely interfere |
| 5. Relations with other people |
| 0   1  2  3  4  5  6  7  8  9   10 |
| Does not interfere                       Completely interfere |
| 6. Sleep |
| 0   1  2  3  4  5  6  7  8  9   10 |
| Does not interfere                       Completely interfere |

| 7. Enjoyment of life | | |
|---|---|---|
| 0           1 2 3 4 5 6 7 8 9           10 | | |
| Does not interfere | | Completely interfere |

In Table 1, there are seven impacts used to evaluate the pain degree of patients. If one patient is interfered by the pain for the sleep, he/she needs to select a score from 0 to 10 to describe the degree of pain. The score *"0"* means he/she is not interfered by the pain while the score *"10"* means he/she is completely influenced.

On the other hand, the treatment cost of lung cancer is high, which will bring economic pressure to patients' families. These pressures might make patients feel guilty about their families and form the stigma (Rush *et al.*, 2005). Because of such stigma, patients might conceal their condition, resulting in inaccurate pain assessment. In this regard, the pain evaluation by nurses should be considered. Hicks *et al.* (2001) gave a useful method to help nurse evaluate the pain degree of patients by their faces. Figure 6 shows the face pain scale given by Hicks *et al.* (2001).

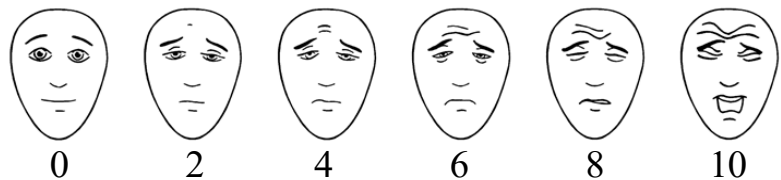

**Figure 6.** The face pain scales (Hicks *et al.*, 2001).

When using the face pain scale, nurses need to observe the patient's face and judge which situation the patient's face belongs to. For instance, if the face of one patient is similar to the second face on the left, the scale of the patient is 2. However, because patients might look different, and the experience of nurses is also different, it might be not accurate to select one of the six faces in Figure 6 to represent the patient.

**4.2. Problem solving**

In Figure 6, the faces of scale 0 and scale 10 are the referenced faces. If the face of one patient is more like the face of scale 0 than the face of scale 10, the pain of the patient is mild. On the contrary, if the face of the patient is similar to the face of scale 10, the pain is serious. In this regard, the similarity of the patient face to the face of scale 0 can be seen as the membership degree while the similarity of the patient face to the face of scale 10 can be seen as the non-membership degree. Then, a CFN consisting of the similarities to the faces of scale 0 and scale 10 is set up.

When giving the CFN in real life, nurses just need to compare the face of patient with the faces of scale 0 and scale 10 to give the similarities of the patient face to the faces of scale 0 and scale 10. For a face of patient, it might not belong to any face in Figure 6, but the similarities of it to the faces of scale 0 and scale 10 can be always given. Therefore, the nurse is able to assess faces with different pain degrees by the CFS. When the experience of nurses is not rich, there might be the overlaps between membership and non-membership degrees which can be modeled by the CFS. In a nutshell, the CFS can be used to express the pain assessments of patients given by nurses.

Suppose one patient uses Table 1 to evaluate the pain, and the normalized average score is 29/70. Then, a nurse assesses the patient and finds that the similarities of the patient face to the faces of scales 0 and 10 are 0.4 and 0.7, respectively. Then, a CFN is set up as $f =< 0.4, 0.7, j >$. Without loss of generality, let $p = 2$ and $\lambda = 0.5$. By Eqs. (4-7), the combined-distance-based score of $f$ is

$$s = \frac{\bar{d}_c(f)}{\bar{d}_c(f) + \tilde{d}_c(f)} \quad \text{where} \quad \tilde{d}_c(f) = 0.5 \times \tilde{d}_m(f) + 0.5 \times \tilde{d}_h(f), \quad \bar{d}_c(f) = 0.5 \times \bar{d}_m(f) + 0.5 \times \bar{d}_h(f),$$

$$\bar{d}_m(f) = \sqrt{(0.4-j)^2 + (0.7-j-1)^2 + j^2 + h^2}, \quad \tilde{d}_m(f) = \sqrt{(0.4-j-1)^2 + (0.7-j)^2 + j^2 + h^2}$$

$\bar{d}_h(f) = \max\{|0.4-j-1|, |0.7-j|\}$, and $\tilde{d}_h(f) = \max\{|0.4-j|, |0.7-j-1|\}$. The higher $s$ is, the less the patient face is similar to the face of scale 10. Hence, $1-s$ shows the pain score of the patient.

To respect the patient's subjective feelings and avoid mistakes caused by the stigma, the gap between the nurse's evaluation $1-s$ and the patient's subjective evaluation 29/70 should be as small as possible. Based on this idea, we establish the linear programming to get $j$ in Programming 1.

**Programming 1.**

$$\min = (41/70 - s)^2$$

$$s.t. \begin{cases} s = \dfrac{\bar{d}_c(f)}{\bar{d}_c(f) + \tilde{d}_c(f)} \\ \tilde{d}_c(f) = 0.5 \times \tilde{d}_m(f) + 0.5 \times \tilde{d}_h(f) \\ \bar{d}_c(f) = 0.5 \times \bar{d}_m(f) + 0.5 \times \bar{d}_h(f) \\ \bar{d}_m(f) = \sqrt{(0.4-j)^2 + (0.7-j-1)^2 + j^2 + h^2} \\ \tilde{d}_m(f) = \sqrt{(0.4-j-1)^2 + (0.7-j)^2 + j^2 + h^2} \\ \bar{d}_h(f) = \max\{|0.4-j-1|, |0.7-j|\} \\ \tilde{d}_h(f) = \max\{|0.4-j|, |0.7-j-1|\} \\ j \in [0.1, 0.4] \end{cases}$$

By solving Programming 1, we can get $j = 0.4$ and $1-s-29/70 = 0.2803$. The result 0.2803

shows two possible situations: 1) The patient conceals the pain, and 2) The nurses overestimate pain. $j$ can represent the confused level of the nurse when he/she assesses. Because $j = 0.4$ is a quite high confused level, the experience of nurse might be not rich, resulting in the second situation. It is suggested to have another nurse to assess the pain of patient. If it is difficult to have another nurse, to avoid the first situation where the patient conceals the pain, the pain score assessed by the nurse $1-s = 0.6946$ is suggested as the final score.

**4.3. Sensitivity analysis**

To demonstrate the reliability of the above results, this section does a simulation to show the changes of $j$ and the gap $41/70 - s$ when $p = 1, 2, 3, 4, 5, 6, 7, 8, 9, 10$ and $\lambda \in [0,1]$. It is found that $j$ is always 0.4. The changes of $41/70 - s$ are shown in Figure 7.

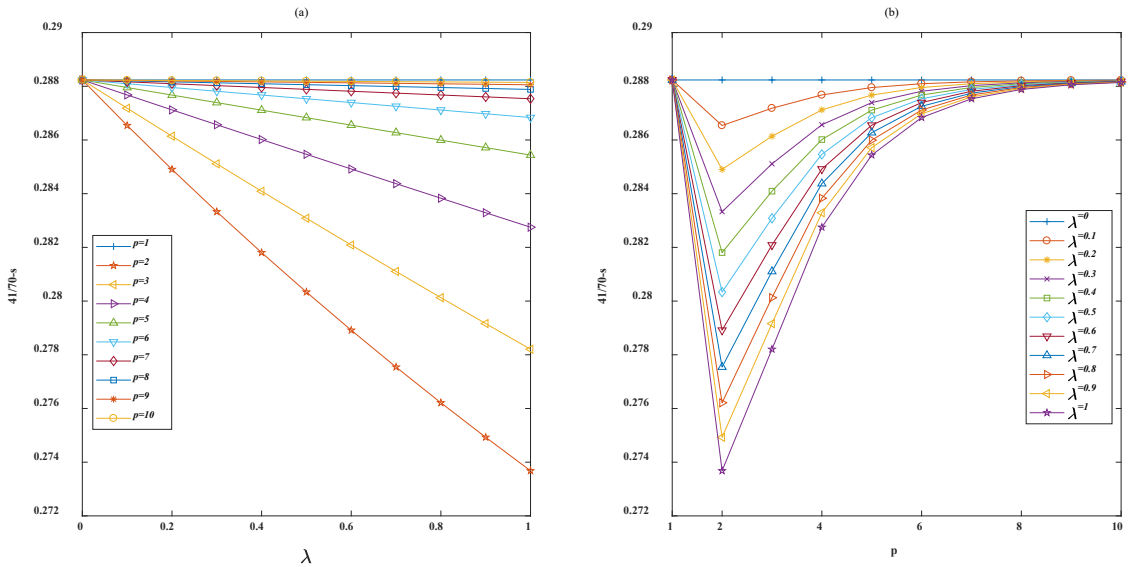

**Figure 7.** The changes of $41/70 - s$ when $p = 1, 2, 3, 4, 5, 6, 7, 8, 9, 10$ and $\lambda \in [0,1]$.

In Figure 7(a), it is easy to find that when $\lambda$ increases with the fixed $p$, the gap $41/70 - s$ descends, but the largest difference between the gaps having $\lambda = 0$ and $\lambda = 1$ is only 0.0045 which is small. Hence, $\lambda$ does not significantly affect $41/70 - s$ in the case study. Known from Figure 7(b), when $p$ ($p \geq 2$) increases with fixed $\lambda$, the gap $41/70 - s$ rises. When $p$ is close to 10, the differences of gaps $41/70 - s$ having various values of $\lambda$ are small. It is clear that when $p = 2$, the smallest gap can be acquired. Therefore, in the lung cancer pain evaluation issue, it is suggested to use $p = 2$ in Programming 1.

## 4.4. Comparison analysis

This section uses Eq. (1), the Minkowski distance proposed by Jiang and Liao (2020), to replace the combined distance in Programming 1. The results are illustrated in Figure 8.

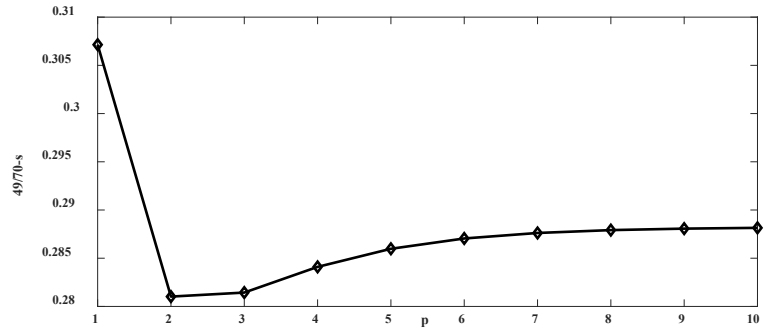

**Figure 8.** The changes of $41/70-s$ when using Eq. (1) with $p=1,2,3,4,5,6,7,8,9,10$.

Because the hesitancy degrees are not considered in Eq. (1), the Minkowski distance proposed by Jiang and Liao (2020) is smaller than the CF-C distance with $\lambda=1$. Hence, the score calculated by Eq. (1) is smaller than the score calculated by the CF-C distance with $\lambda=1$, thereby causing the values of $41/70-s$ in Figure 8 are larger than that in Figure 7(b) when $p$ is the same. From Figures 8 and 7(b), it is easy to find that the difference between the largest value and the smallest value of $41/70-s$ calculated by the Minkowski distance (Jiang & Liao, 2020) is larger than the difference between the largest value and the smallest value of $41/70-s$ calculated by the CF-C distance. This indicates that for different values of $p$, the undulate of the results of the score using the Minkowski distance (Jiang & Liao, 2020) is larger than the undulate of the results of the score using the CF-C distance. Therefore, the CF-C has higher reliability than the Minkowski distance (Jiang & Liao, 2020).

In a nutshell, the comparison analysis shows that Programming 1 based on the CF-C distance can reliably solve the lung cancer pain evaluation issue. Meanwhile, the sensitivity analysis gives an insight that $p=2$ is recommended to be applied in Programming 1 when evaluating the lung cancer pain of patients.

## 5. Conclusions

This paper proposed the CF-IM distance, the CF-H distance, and the CF-C distance. It was found that the anti-perturbation ability of the CF-H distance is higher than that of the CF-IM distance, but the information utilization of the CF-IM distance is higher than that of the CF-H distance. The CF-C distance balanced the CF-IM distance and the CF-H distance by a linear combination of them. When the balance parameter ascended, the anti-perturbation ability of the CF-C distance descended while the information utilization of it increased. The anti-perturbation abilities of the CF-IM distance, CF-H distance, and CF-C distance tended to be similar when $p$ rose. The combined-distance-based score function was introduced and applied to solve the lung cancer pain evaluation problem. The advantages and reliability of the proposed method was verified by the sensitivity and comparison analyses. In addition, it was found that $p = 2$ is suggested to be used when the combined-distance-based score function was employed in evaluating the lung cancer pain of patients.

In the future, the similarity measures and other distance measure of CFS can be studied. The CF-C distance and combined-distance-based score function of CFS can be considered to be employed in other real-life issues.

## Author Contributions

Methodology, Lisheng Jiang; writing-original draft preparation, Tianyu Zhang and Shiyu Yan; writing-review and editing, Ran Fang.

## Funding

The work was supported by the Project of Sichuan System Science and Enterprise Development Research Center (Xq24B03), and the Energy and Environment Carbon Neutrality Innovation Research Center (YB03202408).

## Data Availability Statement

The original contributions presented in this study are included in the article/supplementary material. Further inquiries can be directed to the corresponding author.

## Conflicts of Interest

We confirm that there are no known conflicts of interest associated with this publication and there has been no significant financial support for this work that could have influenced its outcome.